\input amssym.tex

\input pictex.tex

\magnification=1200
\overfullrule =0pt


\newcount\sectnumber%
\def\sectn#1 {\advance \sectnumber by 1{\par\bigskip\noindent{\bfx
\the\sectnumber\ #1}}\par}%

\newcount\subsbnumber%
\def\bsubsn#1 {\advance \subsbnumber by 1{\par\noindent{\cmsc\the\sectnumber.\the\subsbnumber\ 
#1}\quad}}%

\newcount\notenumber

\def\noter#1{\advance\notenumber by 1
$^{\the\notenumber}$\footnote{}{$^{\the\notenumber}$\sevenrm\ #1}}

\newcount\subssnumber%
\def\subssn#1 {\advance \subssnumber by 1{\par\medskip\noindent{\cmsc\the\sectnumber.\the\subssnumber\  #1}\quad}}%

\def\nr#1{\hfill\break {\ninesc#1}\ }%
\def \BR{\Bbb R}
\def\qed{{$\square$}\qquad}
\def\inv#1{#1^{-1}}
\def\rn{\romannumeral}
\font\ninesc=cmcsc9
\def\QED{\hfill{$\square$\smallskip}}%
\outer\def\proclaimsc#1#2{\par\medbreak\noindent{\ninesc #1 \enspace}{\sl
#2}\par\ifdim\lastskip<\medskipamount \removelastskip\penalty55\medskip\fi}%
\def\tup#1 {$#1^{\rm tuple}$ }
\def\tups#1 {$#1^{\rm tuples}$ }
\def\tupsf#1 { $\scriptstyle#1 ^{tuples}$ }
\def\centreline{\centerline}%
\def\frac#1#2{{#1\over #2}}
\def\dfrac{\displaystyle\frac}

\def\recip#1{\frac{1}{#1}}
\def\drecip#1{\dfrac{1}{#1}}
\font\bfxx=cmbx12 at 15pt
\font\sevenit=cmti7
\font\bfx=cmbx12
\font\bfnn=cmbx8
\font\teneufm=eufm10

\font\cmsc=cmcsc10 
\def\eufm#1{{\hbox{\teneufm #1}}}
\def\ref#1{[{\sl #1\/}]}
\def\refb#1{[{\it #1\/}]}%
\def\refbf#1{[{\sevenit #1\/}]}%
\def\GA{{\rm (\!GA\!)}}
\def\func#1#2#3{$#1\colon #2\mapsto #3$}
\def\J#1{$(J_{#1})$}
\def\JJ#1{$(\sim J_{#1})$}

\centreline{\bfxx Accentuate the Negative}\vskip1truecm
 \centreline{\it Dedicated to Professor Pe\v cari\'c on the occasion of his 60th birthday.}
\vskip1truecm
\centreline{\sevenrm {\bfnn Abstract}: A survey of mean inequalities with real weights is given}
\vskip1truecm
\sectn{\ Introduction}
We will be concerned with means that are functions of \tups{n} of real\footnote{}{\sevenrm{\sevenit Mathematics subject classification (2000)}: 26D15}\  numbers with which are associated some positive weights \footnote{}{\sevenrm{\sevenit Keywords}: Convex functions, Jensen inequality, Jensen-Steffensen inequality, reverse inequality}\hskip -.2em,\hskip .2em a typical example being the geometric-arithmetic mean inequality:
$$
\root {W_n\hskip-.3em} \of{a_1^{w_1}\ldots a_n^{w_n}} \le\frac{w_1a_1+\cdots +w_na_n}{W_n},\eqno{\GA}
$$ 
where the weights $w_1, \ldots w_n$ and the variables, $a_1, \ldots, a_n$, are positive numbers and $W_n= w_1+\cdots+w_n$\noter{This notation will be used throughout; given real numbers $ \scriptstyle q_1, q_2,\ldots q_n$ then $ \scriptstyle  Q_k = \sum_{i=1}^k q_i,\,1\le k\le n$. Also we write $\scriptstyle \tilde Q_k =$\hfill\break\vskip-.45truecm $\scriptstyle Q_n-Q_{k-1} = \sum_{i= k}^na_i,\,1\le k\le n. $}. There is no real reason for excluding zero values for the weights  except that if for instance $w_n=0$  this effectively  means that we are stating  or discussing  the inequality for a smaller value of $n$; equivalently $\GA$ then states the inequality for all values of $k, 1\le k\le n$. A similar remark can be made about assuming all the variables are distinct.

However  it is usual  not to allow negative weights even though there is a very good and useful theory that covers this possibility.  Classically the first person to  study this  in detail was   Steffensen early in  the twentieth century.  More recenty very significant contributions have been made by  Pe\v cari\'c and his colleagues. The case of real weights  has been of interest to Pe\v cari\'c throughout his career from his  student days  up to the present. However the results are not generally known  and  this paper is  an attempt to remedy this neglect. 

Since almost  all the inequalities between means are particular cases of the Jensen inequality for convex functions\noter{Thus $\scriptstyle\GA$ is just a property of the convexity of the  function $ \scriptstyle f(x) = -\log x$, or the convexity of $\scriptstyle g(x) = e^x$; \refbf{B03 p.92, BB\hfill\break\vskip-.55truecm  pp.6-7}.}the paper will concentrate on this result. Applications to particular means will then follow using  the lines of the original application of Jensens's inequality.
\sectn{\ Convex Functions} The definitions and properties of convex functions are well known and will not be given in detail here.  However the basic inequality of Jensen is equivalent to the definition of convexity and so in this section we will give details that are necessary for later discussion. 

Perhaps the simplest analytic  definition of a convex function is: let $I$ be an open interval $I\subseteq \Bbb R$\noter{This meaning for $ \scriptstyle I$ will be used throughout the paper.} then \func{f}{I}{\Bbb R} is convex  if  $\forall x,y\in I$ then \func{D}{]0,1[}{\Bbb R}  
 is non-positive,
$$
D(t)= D_2(t)= f\bigl((1-t)x + ty\bigr)-\bigl((1-t)f(x)+tf(y)\bigr)\le 0.\eqno(1)
$$
It should be noted that if $x, y \in I$ then so is $\overline  x=(1-t)x + ty,\; \forall\, t,\, 0< t<1$, so all the terms on the right-hand side are defined\noter{More precisely $ \scriptstyle \overline x \in I_0= [\min\{x,y\},\,  \max\{x, y\}]$.}. Further note that $D$ is defined for all $t\in 
\Bbb R$ and use will be made of this in later discussions.

An alternative but equivalent definition is:    $\forall\; z\in I$ there is an affine function \func{S_z}{\Bbb R}{\Bbb R} such that:
$$
S_z(z)= f(z)\quad \hbox{and}\quad S_z(x)= f(z) + \lambda(x-z) \le f(x)\;\; \forall  x\in I.
$$
See \refb{B03 p. 27, HLP pp.70--75, 94--96,  PPT p.5, RV p.12}.

The geometric interpretations of these definitions are immediate from Figures 1 and 2. 
\bigskip
\beginpicture
\setcoordinatesystem units <1.3cm,  1.3cm>
\setlinear \plot  0 0 2 2 /
\setlinear \plot  4.2 -1 7.5 -1 /
\setlinear \plot  5.5 -1 5.5 1.6 /
\setquadratic \plot   0 .7 1 .5 1.7 2.1 /
\setquadratic \plot   4.5 1   6 -1.4  6.8 .1 /
\setlinear \plot  .4 .4   .4 -1  /
\setlinear \plot 1.56 1.6   1.56 -1  /
\setlinear \plot  -1 1.5   1.59 1.55   /
\setlinear \plot  -1 .43   .39 .43   /
\setlinear \plot  -1 -1   2 -1   /
\setlinear \plot  -1 -1   -1 1.9   /
\setlinear \plot  .79 -1  .79 .8   /
\setlinear \plot  1.3 -1  1.3 .87   /
\setlinear \plot  1 .4  1.9  1.78 /
\put {$\scriptstyle S_z$} at 2. 1.6
\put {$\scriptstyle z$} at 1.3 -1.1
\put {$\scriptstyle \overline x$} at .79 -1.1
\put {$\scriptstyle x$} at .39 -1.1
\put {$\scriptstyle y$} at 1.56 -1.1 
\put {$\scriptstyle Graph\ of \ D$} at 6.15 0.6
\put {$\scriptstyle 0$} at 5.49 -1.15
\put {$\scriptstyle 1$} at 6.57 -1.13
\put {$\scriptstyle Graph\ of \ a\ Convex Function$} at .24 2.3
\put {$\scriptstyle f(x)$} at -1.3 .38
\put {$\scriptstyle f(y)$} at -1.3 1.5
\put {$\scriptstyle Q$} at .7 .9
\put {$\scriptstyle P$} at .88 .27
\put {$\scriptstyle P=(\overline x, f(\overline x)),$} at -1 -1.5
\put {$\scriptstyle Q=(\overline x,f((1-t)x+ty)) ),$} at 1 -1.5

\endpicture
\bigskip
Figure 1\hskip 7cm Figure 2  
\bigskip
Use will be made of the following properties of convex functions.

(C1)  The first divided difference $[x,y; f]= \bigl(f(x) - f(y)\bigr)\big/(x-y),\, x,y\in I,\, x\ne y,$ is increasing in both variables; \refb{B03 p.26, PPT p.2, RV p.6}.

(C2) If $x,y,z,u,v\in I$ and $x\le y\le z\le u \le v$ and if $S_z(t) = f(z) + \lambda (t-z)$ then: 
$$
f(y)-f(x)\le \lambda(y-x),\qquad f(v)-f(u)\ge \lambda(v-u).
$$
See \ref{P84b}.

(C3) A  function convex on $I$ is continuous; [RV p.4].\hskip-.1em \noter{But not necessarily differentiable; consider $ \scriptstyle f(x) =\vert x\vert$.} 
 
(C4) The Hardy-Littlewood-P\'olya--Karamata-Fuchs  majorization theorem, or just HLPKF, \refb{B03 pp.23, 24, 30, BB pp.30--32, HLP pp.88-91, PPT pp.319--320, MO pp. 64--67}: if ${\bf a}=(a_1, \ldots, a_n)$\noter{\sevenrm This notation for \tupsf{n} or sequences, will be used throughout.}, ${\bf b}=(b_1, \ldots, b_n) $  are decreasing  \tups{n} with entries in the domain of a convex function $f$ and ${\bf w}=(w_1, \ldots, w_n)$ a real \tup{n}  and if:
$$
\sum_{i=1}^k w_1a_i\le \sum_{i=1}^k w_1b_i,\, 1\le k<n\quad \hbox{and}\quad \sum_{i=1}^n w_1a_i=\sum_{i=1}^n w_1b_i
$$
 then:
$$
\sum_{i=1}^n w_1f(a_i)\le \sum_{i=1}^n w_1f(b_i).
$$
 The first two are rather elementary and have obvious geometric interpretations  but the last two are more sophisticated.

Jensen's  inequality is an easy deduction from the definition of convexity and in a variety of forms is given in the following theorem.
\proclaimsc{Theorem 1}{ Let $n\in \Bbb N, n\ge 2, I$ an  interval, \func{f}{I}{\Bbb R} convex then:
$$
\eqalignno
{ 
(a)\;\forall x_i\in I, 1\le i\le n,&\;\hbox{and}\;\forall t_i,1\le i\le n, \;\hbox{such that}\; 0< t_i<1, \; 1\le i\le n,\cr  
\hbox{and}\; t_1= 1-\sum_{2}^nt_i&\;\hbox{we have:}\cr
D(t_2, \ldots t_n) = &D_n(t_2, \ldots t_n)= 
f\Bigl( \sum_{i=1}^nt_ix_i \Bigr)-\sum_{i=1}^nt_if(x_i) 
\le 0;\cr
(b)\; \forall a_i \in I,\, 1\le i\le n,& \; \hbox{and}\; \forall\; \hbox{positive weights}\; w_i,1\le i
\le n,\phantom{aaaaaaaaqqqqqqqqqq}\cr 
 &f \Biggl(\recip{W_n}\sum_{i=1}^nw_ia_i \Biggr)\le\recip{W_n} \sum_{i=1}^nw_if(a_i);& \Biggr\rbrace(J_n)\cr
(c)\;\forall\; a_i\in I,\, 1
\le i\le n,& \; \hbox{and}\;\hbox{positive weights}\; p_i,\, 1\le i\le n,\;\hbox{with}\; P_n =1,\phantom{aaaaaaaaq}\cr  
&f\Bigl(\sum_{i=1}^np_ia_i\Bigr)\le \sum_{i=1}^np_if(a_i).\cr
}
$$
}\par
\qed {\it Proof (\rn1)} The most well known proof is by induction, the case $n=2$, \J{2}, being just (1), a definition of convexity; \refb{B03 p.31, PPT pp. 43--44}, \ref{P90}. 

{\it Proof (\rn2)} Another proof can be based on the support line definition above; \refb{RV pp.189--190}, \ref{P84b}.

{\it Proof (\rn3)} A geometric  proof can be given as follows.  

First note, using (1),  that the set bounded by the chord joining $\bigl(x, f(x)\bigr)$ to $\bigl(y, f(y)\bigr)$   and the graph of $f$ joining the same points is a convex set. Then by induction show that the point $ (\overline a, \overline \alpha)$\noter{This point is just the weighted centroid of the points $ \scriptstyle (a_i, f(a_i)),\, 1\le i\le n$, that lie on the graph of $ \scriptstyle f$.}, $ \overline a = \sum_{i=1}^np_ia_i, \overline \alpha = \sum_{i=1}^np_if(a_i)$, lies inside this set and so $\overline \alpha\ge f(\overline  a)$ which is just \J{n}.
\QED
\smallskip

We now turn to the main interest of this paper. What happens if we allow negative weights in \J{n}?
\sectn{\ The Case of Two Variables}
The inequality \J{2} is just $D(t)\le 0,  0< t< 1,$ and 
it is immediate from Figures 1 and 2 that if either $t< 0$ or $1-t< 0$, equivalently $t>1$, then $D(t)\ge 0$, that is the reverse inequality\noter{\sevenrm The naming of reverse inequalities varies; sometimes the term inverse is used and  sometimes converse but reverse seems to be the\hfill\break\vskip-.55truecm best usage,} holds.  Formally we have the following result where the last of the  notations in Theorem 1 is used, \refb{B03 p.33}, \ref{LP}.
\proclaimsc{Theorem 2}{If $f$ is convex on the  interval I and either $p_1< 0$ or $p_2< 0$ then forall  $a_1, a_2$ with $\overline a= p_1a_1 + p_2a_2\in I$
$$
f(p_1a_1+ p_2a_2)\ge p_1f(a_1)+ p_2f(a_2). \eqno(\sim J_2)
$$} 
\qed {\it  Proof (\rn1)}   It is an easy exercise to use the second definition of convexity to prove that the function $D$ is convex on $\Bbb R$.  Hence since $D(0)= D(1)= 0$ we must have that $D(t)\le 0, 0<t<1$ and $D(t)\ge 0, t< 0,\,  t> 1$, as shown in Figure 2.

{\it Proof (\rn2)} Assume that $p_2< 0$ then: 
$$
\eqalign
{
a_1 =& \frac{\overline a -p_2a_2}{p_1}= \frac{\overline a -p_2a_2}{1-p_2}.\cr
\hbox{So, using $(J_2)$,}\cr
f(a_1) =& f\Bigl(\frac{\overline a -p_2a_2}{1-p_2}\Bigl)\cr
\le&\frac{f(\overline a) -p_2f(a_2)}{1-p_2}=\frac{f(\overline a) -p_2f(a_2)}{p_1}.\cr
}
$$
Rewriting the last line gives $(\sim J_2)$.

{\it Proof  (\rn3)} Let us assume that $t< 0$ and without loss of generality that $a_1\le a_2$.   

Then $a_1$ lies between $ \overline a$ and $a_2$ and $ a_1= (\overline  a - t a_2)\big/(1-t)$. Now let $S= S_{a_1}$ then
$$
\eqalign
{
f(a_1) = &S(a_1) =S\Bigl(\frac{\overline a -ta_2}{1-t}\Bigr)=\frac{S( \overline a)-tS(a_2)}{1-t}\cr
\le &\frac{f( \overline a) -tf(a_2)}{1-t},\cr
}
$$
which on rewriting gives \JJ{2}.\QED

Note that the condition $ \overline a\in I$ is necessary as the $ \overline a\notin I_0$ and so may very well not lie in the domain of $f$ \noter{Clearly if $\scriptstyle I= \Bbb R$ the condition can be omitted.}.

In the case of two variables \ the situation is completely determined: either the weights  are positive when we have Jensen's inequality or one is negative when we have the reverse of inequality\noter{The geometric-arithmetic mean inequality case of this result was the motivation for one of Pe\v cari\'c's more interesting collaborations.}.  Matters are not so simple when we have three or more variables.

In other terms: for all $x,y \in I$ with $ \overline x
\in I$ the sets $D_+=\{t;\, t\in \Bbb R \land D(t)>0 \}, D_-=\{t;\, t\in \Bbb R \land D(t)<0 \}, D_0=\{t;\, t\in \Bbb R \land D(t)=0 \}$ partition $\Bbb R$ and do not depend on $x$ or $y$.

This very simple result has been given this much attention as the ideas and  methods of proof are used in the more complicated cases of more than two variables.

\sectn{The Three Variable Case}  
This case is very different to the two variable situation discussed above but has its own peculiarities; in addtion it introduces ideas needed for the general case. The function $D$ can now be written:
$$
D(s,t) = D_3(s,t) =f\bigl((1-s-t)x+sy+tx\bigr)-\bigl((1-s-t)f(x) + sf(y) +tf(z)\bigr).
$$ 
Clearly $D_3$ partitions $\Bbb R^2$ into three sets\noter{Using the notation of the previous section.}: the closed convex $0$-level curve $D_0$, the open convex set $D_-$, that is the interior of this curve and where \J{3} holds, and the unbounded exterior of the this curve, $D_+$, where \JJ{3} holds.  However unlike the two variable cae these sets depend on the other variables $x,y,z$ as we will now see.

 The set where Jensen's inequality, $(J_3)$, holds for all $x,y,z \in I$, is the   triangle $T$  where the above weights are positive 
$$
T=\{(s,t);\,  0< s<1,\; 0< t< 1, 0< s+t < 1\};
$$  see Figure 3. 

\beginpicture
\setcoordinatesystem  units <1.3cm, 1.3cm> point at 0 0
\setlinear \plot 5 -4 5 1  /
\setlinear \plot  3 -2  8 -2 /
\setlinear \plot  7 -3  4 0 / 
\setlinear \plot  7 -4  3.5 -.5 / 
\setlinear \plot  3 -1  7 -1 /
\setlinear \plot  6 -4  6 -0 /
\put {$\scriptstyle S_2$} at 5.8 -1.3
\put{\rm Figure 3} at 9.5 -1.3
\put {$\scriptstyle S_1$} at 4.8 -1.3
\put {$\scriptstyle S_3$} at 5.8 -2.3
\put {$\scriptstyle T$} at 5.3 -1.75
\put {$\scriptstyle T_1$} at 4.5  .3
\put {$\scriptstyle T_2$} at 7  -2.4
\put {$\scriptstyle T_3$} at 4  -2.9
\put {$\rightarrow$} at 8 -2.02
\put {$\uparrow$} at 5 1
\put {$ \scriptstyle s$} at 8.1 -2.2
\put {$ \scriptstyle t$} at 5.2 1
\put {$ \scriptstyle (0,0)$} at  4.7 -2.19
\put {$ \scriptstyle (1,0)$} at  6.3 -1.88
\put {$ \scriptstyle (0,1)$} at  5.3 -.8
\put {$\scriptstyle s+t =1$} at   7.5 -2.9
\put {$\scriptstyle s+t =0$} at   7.5 -3.9
\put {$\scriptstyle s =1$} at   6  .1
\put {$\scriptstyle t =1$} at   2.7  -1
\endpicture 
\bigskip
On the sides of this triangle one of the weights is zero and so we have cases of the two variable situation, as we noted above, and as a result  by\J{2} $D_3\le 0$ on the sides of  $T$.  Hence by $(C_3)$ $D_3$ must be negative on a set larger than the triangle,  that is $T\subset D_0$; note that the vertices of $T$ lie on $D_0$. In any case for some choices of $x,y,z\in I$ \J{3}  holds with negative weights and the question is  whether there is a larger set than $T$ on which $(J_3)$  holds for  a large  choice of variables, or for variables satisfying some simple condition: \refb{B03 pp. 39--41}, \ref{B98}. 

Let us look at what happen when negative weights are allowed.  

The next result, due to Pe\v cari\'c, \ref{P81a, P84a, VP}, resolves the case when there is a maximum number of negative weights:\noter{Clearly three negative weights is the same  as three positive weights.}\refb{B03 p.43,  PPT p.83} .    

\proclaimsc{Theorem 3}{If $f$ is convex on the interval I and only one of $p_1, p_2, p_3$ is positive and if $a_1, a_2, a_3,  \overline a= p_1a_1 + p_2a_2+ p_3a_3\in I$ then
$$
f(p_1a_1+ p_2a_2+p_3a_3)\ge p_1f(a_1)+ p_2f(a_2)+ p_3f(a_3). \eqno(\sim J_3)
$$} 
\qed  {\it  Proof (\rn1)} If we consider $D(s,t),\, (s,t)\in \Bbb R^2$, and assume $f$ is differentiable then it can easily be shown that $D$ has no stationary points in its domain.  An immediate conclusion is that \J{3} must hold in the triangle $T$ since the maximum and minimum of $D$ must occur on the boundary and it is non-positive there by $(J_2)$.  The domains where  two of the weights are negative are the three unbounded triangles  $T_1, T_2, T_3$ of Figure 3.  By Theorem 2 $D$ is non-negative  on the boundaries of these triangles  and so it would be reasonable to conclude that $D$ is non-negative on these triangles giving a proof of \JJ{3}.  This proof is not quite complete as these are unbounded regions and this simple argument does not work. Let us   look at the second proof of Theorem 2.

{\it  Proof (\rn2)}  Assume without loss in generality that $p_1>0, p_2< 0, p_3<0$ then 
$$
\eqalign
{
a_1 =& \frac{\overline a -p_2a_2-p_3a_3}{p_1}= \frac{\overline a -p_2a_2-p_3a_3}{1-p_2-p_3}.\cr
\hbox{So, using \J{3},}\cr
f(a_1) =& f\Bigl(\frac{\overline a -p_2a_2-p_3a_3}{1-p_2}\Bigl)\cr
\le&\frac{f(\overline a) -p_2f(a_2)-p_3f(a_3)}{1-p_2}=\frac{f(\overline a) -p_2f(a_2)-p_3f(a_3)}{p_1}.\cr
}
$$
Rewriting the last line gives \JJ{3}.\QED

\smallskip

It remains to consider what happens if there is only one negative weight.   In order for  \J{3} to hold we need $ \overline a\in I_0$, and for \JJ{3}to hold  $ \overline a\notin I_0$.   Assume without loss in generality that $a_1<a_2<a_3$ and assume that $p_1<0$ then 
$$
\overline a= p_1a_1 + (p_2 + p_3) \frac{a_2p_2 + a_3p_3}{p_2 + p_3},
$$
The second term on the right of the last term is in the interval $]a_2,a_3[$ and so $ \overline  a$ is to the right of $a_2$ and can lie either in $I_0$ or not depending on the value of the negative $p_1$. Further any condition on $p_1$ to require one or other of these options would obviously depend on thevalues of $a_1,a_2, a_3$.

A similar argument applies if the negative weight is $p_3$.

However in the case of the middle term $a_2$ having a negative weight, $p_2<0$.  Steffensen, [S],  obtained a simple condition on the weights that would assure $\overline a\in I_0$. Consider
$$
\eqalign{
\overline a=& p_3(a_3-a_2) + (p_3+p_2)(a_2-a_1) + a_1\cr
=& p_1(a_1-a_2) + (p_1+p_2)(a_2-a_3) + a_3.\cr
}
$$ 
 If we  assume that $ p_3+p_2>0$ the first expression shows that $\overline a> a_1$ and if we require that $p_1+p_2>0$ the second expression shows that $\overline a< a_3$. That is: with these two conditions on the weights  $\overline a\in I_0$ and  \J{3}  should hold.

The conditions can be put in a simpler form:
$$
0<p_1<0,\quad 0<P_2= p_1 + p_2<1\eqno(S_3)
$$
since $(S_3)$ is easily seem to be equivalent to  
$$
0<p_3<0,\quad 0<\tilde P_2= p_3 + p_2<1\eqno(\tilde S_3)
$$
Later Pe\v cari\'c, [P81a], gave an alternative form of this condition: the negative weight is dominated by both of the positive weights, that is
$$
P_2>0, \quad\tilde P_2>0,\eqno(P_3)
$$
Hence we have the following result of Steffensen, [S]; several proofs are given, in addition to the one sketched above since they extend to give different results when $n>3$,

\proclaimsc{Theorem 4}{If $f$ is convex on the interval I and if $0<p_1<1,0< P_2< 1$  and either $a_1\le a_2\le a_3$ or $a_1\ge a_2\ge a_3$,with $a_1, a_2, a_3\in I$, then
$$
f(p_1a_1+ p_2a_2+p_3a_3)\le p_1f(a_1)+ p_2f(a_2)+ p_3f(a_3). \eqno(J_3)
$$}\par 
\qed{\it  Proof (\rn1)} \refb{B03 p. 39} Without loss in generality assume that $a_1< a_2< a_3$ and  write $\tilde a= P_2a_2 + p_3a_3$; note that $a_1< \overline a< a_3$ and $a_2< \tilde a<a_3$.
Now 
$$
\eqalign
{
p_1f(a_1)+p_2f(a_2)+p_3f(a_3)-f( \overline  a)& =-p_1\bigl(f(a_2)-f(a_1)\bigr)+P_2f(a_2) +p_3f(a_3)-f( \overline a)\cr
&\ge -p_1\bigl(f(a_2)-f(a_1)\bigr) + f(\tilde a)-f(\overline a).\quad\hbox{by $(J_2)$}, \cr
&=p_1(a_2-a_1)\biggl(\frac{f(\tilde a)-f( \overline a)}{\tilde a- \overline a}-
\frac{f(a_2)-f(a_1)}{a_2-a_1}\biggr)\cr
&\ge 0,\quad \hbox{by $(C_1)$}.\cr
}
$$

{\it  Proof (\rn2)} \ref{P81a} In this proof we use the fact  $p_2<0$, this was not used in the first proof. 

Again asuming   $a_1< a_2< a_3$ we have that for some $t,\, 0<t<1$ that $a_2=(1-t)a_1+ta_3$. Then:
$$
\overline a = p_1a_1 + p_2(1-t)a_1+ta_3 + p_3a_3= \bigl(p_1+(1-t)p_2\bigr)a_1+(p_3+tp_2)a_3.
$$
So by $(J_2)$,  noting that the coefficents of the last expression are positive and have sum equal to 1,
$$
\eqalign
{
f(\overline a) = &f\Bigr(\bigl(p_1+(1-t)p_2\bigr)a_1+(p_3+tp_2)a_3\Bigr)\cr
\le&\bigl(p_1+(1-t)p_2\bigr)f(a_1)+(p_3+tp_2)f(a_3)\cr
=&p_1f(a_1) +p_2\bigl((1-t)f(a_1) + tf(a_3)\bigr) + p_3f(a_3)\cr
\le&p_1f(a_1) +p_2f\bigl((1-t)a_1+ta_3\bigr) + p_3f(a_3)\cr
=&p_1f(a_1) +p_2f(a_2) + p_3f(a_3),\cr
}
$$
{\it  Proof (\rn3)} \refb{PPT pp. 57--58}, \ref{P84a}  Assume without loss in generality that $a_1<\overline a<a_2<a_3$
 and  define $ \lambda$ by: 

$S_{\overline a}(x)= f(\overline a)+ \lambda(x- \overline a)$.  Using (C2)  we get: 
$$
\eqalign
{
p_1f(a_1)+p_2f(a_2)+&p_3f(a_3)-f(\overline a)\cr=&p_1\bigl(f(a_1)-f(\overline a)\bigr) +(p_2+p_3)\bigl(-f(\overline a)+f(a_2)\bigr) + p_3\bigl(f(a_3)-f(a_2)\bigr)\cr
&\ge p_1 \lambda(a_1- \overline a)+(p_2+p_3)\lambda (a_2- \overline a) + p_3\lambda(a_3-a_2) = 0.
}
$$

{\it  Proof (\rn4)} \ref{P84b} Without loss in generality assume that $b_1=a_1>b_2=\overline a>b_3=a_2>b_4=a_3$ further define $q_1= p_1, q_2= -1, q_3=p_2, q_4= p_3$; then if $c_i= \overline a, 1\le i \le 4$ :
$$
\displaylines
{
q_1b_1=p_1a_1\ge p_1 \overline a = q_1c_1\cr
q_1b_1 + q_2 b_2=p_1a_1-\overline a\ge q_1c_1 + q_2 c_2\cr
q_1b_1 + q_2 b_2+q_3b_3=p_1a_1+p_2a_2-\overline a\geq_1c_1 + q_2 c_2+ q_3c_3\cr
q_1b_1 + q_2 b_2+q_3b_3 +q_4b_4=0 =q_1c_1 + q_2 c_2+ q_3c_3 + q_4c_4
} 
$$
hence by $(C_4)$, HLPKF:
$$
q_1f(b_1) + q_2 f(b_)+q_f(b_3) +q_4f(b_4)\ge q_1f(c_1) + q_2f(c_2)+ q_3f(c_3) + q_4f(c_4)= 0
$$
which is just $(J_3)$.

{\it  Proof (\rn5)} The Steffensen condition  tells us that  the point $\bigl( \overline a,p_1f(a_1) +p_2f(a_2) + p_3f(a_3)\bigr)$ lies in the convex hull of the points $\bigl(a_i, f(a_i)\bigr),\; a\le i\le 3$, and so lies  in the convex set $\{(x,y); y\ge
f(x)\}$ and this implies $(J_3)$.\QED 

Using the notation in the definition of $D_3$ and assuming that $x<y<z$ and $ s<0$ the condition $(S_3)$ is just: $0\le s+t\le 1, 0\le t\le 1$ and so $(J_3$) holds in the triangle $S_1$ of Figure 3.  Depending on the order of $x,y,z$ and provided the central element has the only negative weight and $(S_3)$ holds then $(J_3)$ will hold in one of $S_1, S_2, S_3$ of Figure 3.

\sectn{The $\bf n$ Variable Case}

 In this section we turn to the general situation and  the notations are those of Theorem 1.

	Let us first consider the extension of Theorem 3. The second proof of Theorem 3 can easily be adapted to the following result of Pe\v cari\'c; \refb{B03 p.43, PPT p. 83}, \ref{VP}.

\proclaimsc{Theorem 5}{If \func{f}{I}{\Bbb R} is convex, $n\in \Bbb N,\, n\ge 2$, $a_i\in I,\,  w_i\in \Bbb R, w_i\ne 0$, $1\le i\le n$, further assume that all the weights are negative except one, $W_n\ne 0$,   and that $\overline a\in I$ then:
$$
f \Biggl(\recip{W_n}\sum_{i=1}^nw_ia_i \Biggr)\ge\recip{W_n}\sum_{i=1}^nw_if(a_i),
$$
or, using an alternative notation, \hfill$\Biggr\rbrace$\JJ{n}
$$
f\Bigl(\sum_{i=1}^np_ia_i\Bigr)\ge \sum_{i=1}^np_if(a_i).
$$
}\par 
\qed  The case $n=2$ is Theorem 2,  and the case n =3 is Theorem 3. 

Assume then $ n\ge 3$ and, without loss in generality, that $p_1>0$ and $p_i<0,\, 2\le i\le n,$ then: 
$$
a_1= \frac{\overline a +\sum_{i=2}^n(-p_i)a_i}{p_1}=\frac{\overline a +\sum_{i=2}^n(-p_i)a_i}{1+\sum_{i=2}^n(-p_i)}.
$$ 

So by \J{n},
$$
\eqalign
{
f(a_1)\le &\recip{1+\sum_{i=2}^n(-p_i)}\Bigl(f( \overline a) +\sum_{i=2}^n(-p_i)f(a_i)\Bigr)\cr
=&\recip{p_1}\Bigl(f( \overline a) +\sum_{i=2}^n(-p_i)f(a_i)\Bigr),\cr
}
$$
which on rewriting is just \JJ{n}.\QED

We now turn to the situation where  \J{n} holds but there are negative weights, the generalization of Theorem 4 due  Steffensen. Note that from Theorem 5 we will need at least two positive weights for \J{n} to hold.

The important conditions put on the weights  of Steffensen and Pe\v cari\'c, $(S_3)$, and $(P_3)$ above,   now differ and are as follows, using the alternative notion of Theorem 5.
\medskip
(S)\qquad $0<P_i<1,\, 1\le i\le n-1$; and of course $P_n=1$.

This implies that  $0<\tilde P_k<1,\,1< k \le n$, and in particular that $0<p_1<1$ and $0<p_n<1$. 
\smallskip
For (P) we introduce the following notation:

\ \phantom{aaaaaaaaa} $I_+=\{i;1\le i\le n\, \land\, p_i >0\}$ and $I_-=\{i;1\le i\le n\, \land\, p_i <0\}$;
\smallskip
obviously $I_+\cap I_=\emptyset$ and $I_+\cup I_-=\{1,2,\ldots, n\}$.
\medskip
(P)\qquad $p_1, p_n\in I_+$  and $\;\forall\, i\in I_+,\;  p_i+\sum_{j\in I_-}p_j>0$.

It is easy to see that (P) implies (S). Further we have the following simple result, \refb{B03 p.38, PPT pp.37--38}.

\proclaimsc{Lemma 6} {If $\bf a$ is monotonic and  (S) holds then $\overline a\in I_0$.}\par
\qed Assume without loss in generality that the \tup{n} is increasing when since
$$
\overline a=\sum_{i=1}^np_ia_i= a_n + \sum_{i=1}^{n-1}P_i(a_i-a_{i+1})
$$

the result follows by (S).\QED

All the proofs of Theorem 4 can be extended to give a proof of the general case.

\proclaimsc{Theorem 7}{Let $n\in \Bbb N, n\ge 3, I\subseteq \overline \Bbb R$ an interval, $f\colon I\mapsto \Bbb R$ convex then  for all  monotonic \tups{n} with terms in $I$ \J{n} holds  for all non-zero real weights satisfying (S).}\par
\qed 
{\it Proof (\rn1)} The standard proof is by induction starting with then case $n=3$,  Theorem 4; see \refb{B03 pp.37--39}.

{\it  Proof (\rn2)} This proof, due to Pe\v cari\'c, [P81a],  needs only the weaker condition (P) and we also assume without loss in generality that the \tup{n} is increasing and distinct.    

 If $i\in I_-$ then $a_1<a_i<a_n$ and hence for some $t_i,\, 0<t_i<1$, $a_i = (1-t_i)a_1 + t_ia_n$ and so
$$
\eqalign
{
\overline a=&\sum_{i\in I_+}p_ia_i+\sum_{i\in I_-}p_i\bigl((1-t_i)a_1 + t_ia_n\bigr)\cr
=&\Bigl(p_1 +\sum_{i\in I_-}p_i\bigl((1-t_i)\Bigr)a_1 +\sum_{i\in I_+\setminus\{1,n\}}p_ia_i +\Bigl(p_n +\sum_{i\in I_-}p_it_i\Bigr)a_n\cr
}
$$
Note that the sum of the weights in this last expression is 1 and that by (P) they are all positive.  Hence by Jensen's inequality
$$
\eqalign
{
f(\overline a)\le& \Bigl(p_1 +\sum_{I_-}p_i\bigl((1-t_i)\Bigr)f(a_1) +\sum_{I_+\setminus\{1,n\}}p_if(a_i) +\Bigl(p_n +\sum_{I_-}p_it_i\Bigr)f(a_n)\cr
=&\sum_{i\in I_+}p_if(a_i)+\sum_{i\in I_-}p_i\bigl((1-t_i)f(a_1)+t_if(a_n)\bigr)\cr
\le &\sum_{I_+}p_if(a_i)+\sum_{I_-}p_if(a_i),\hbox{by $(J_2)$ and the negativity of the $p_i$ in the last sum}\cr
=&\sum_{i=1}^n p_if(a_i).
}
$$
which is  $(J_n)$.

{\it  Proof (\rn3)} \refb{PPT pp.57--58} Assuming without loss in generality that the \tup{n}  is distinct  and  decreasing we have from  Lemma 6 that $a_n< \overline a< a_1$ and assume that $a_k\ge \overline a\ge a_{k+1}.\, 1\le k< n$. As in Theorem 4 define $ \lambda$ by: 

$S_{\overline a}(x)= f(\overline a)+ \lambda(x- \overline a)$; then by (C2):
$$
\overline a\le u\le v\Longrightarrow f(v) -f(u)\ge \lambda(v-u);\quad u\le v\le \overline a\Longrightarrow f(v) -f(u)\le \lambda(v-u).
$$ 
By Lemma 6 $a_1\ge \overline a\ge a_n$ so assume first that $a_1> \overline a> a_n$, and that $a_{k+1}\le  \overline a\le a_k$  for some $k\, 1\le k\le n-1$. Then
$$
\eqalign
{
f(\overline a)-&\sum_{i=1}^np_if(a_i)=f(\overline a)-\sum_{i=1}^kp_if(a_i)-\sum_{i=1}^{k+1}p_if(a_i)\cr
=&f(\overline a)-\sum_{i=1}^{k-1}P_i\bigl(f(a_i)-f(a_{i-1})\bigr) -P_kf(a_k)-\sum_{i=k}^{n-1}\tilde P_i\bigl(f(a_{i+1}-f(a_i)\bigr) -\tilde P_{k+1}f(a_k)\cr
=&\sum_{i=1}^{k-1}P_i\bigl(f(a_{i-1})-f(a_i)\bigr)+P_k\bigl(f(\overline a)-f(a_k)\bigr)+\tilde P_{k+1}\bigl(f(\overline a)-f(a_k)\bigr)+\sum_{i=k}^{n-1}\tilde P_i\bigl(f(a_i)-f(a_{i+1})\bigr)\cr
\ge&\sum_{i=1}^{k-1}\lambda P_i(a_{i-1}-a_i)+ \lambda P_k( \overline a -a_k)+ \lambda \tilde P_{k+1}( \overline a -a_{k+1})+\sum_{i=k}^{n-1}\lambda\tilde P_i(a_i-a_{i+1})\cr
=& \lambda\bigl(\overline a- \sum_{i=1}^np_ia_i\bigr)=0.\cr
}
$$
{\it  Proof (\rn4 )} \ref{P81a}
Using the notations and assumptions of the previous proof  define  the three \tups{(n+1)}  $ x_1,\ldots, x_{n+1}, y_1, \ldots, y_{n+1}, q_1, \ldots q_{n+1}$:
$$
\eqalign
{
x_i = a_i,&\; q_i = p_i,\quad 1\le i\le k;\cr
x_{k+1} = \overline a,&\; q_{k+1}=-1;\cr
x_i=a_{i-1}, &\; q_i = p_{i-1},\quad k+2\le i\le n+1;\cr
y_i= \overline a,&\quad 1\le i\le n+1.
}
$$
Simple calculations show that: $Q_j = P_j,\; 1\le j\le k, = \tilde P_{j-1},\; k+1\le j\le n, = 0, j=n+1 $; and, 
$$
\displaylines
{
\sum_{i=1}^jq_iy_i=P_j \overline a,\quad 1\le j\le k,\cr
\phantom{aaaaaaaaaaaaaaa}=\tilde P_{j-1}\overline a,\quad k+1\le j\le n,\cr
\phantom{aaaaaaa}=0,\quad j=n+1.\cr
\sum_{i=1}^jq_ix_i=\sum_{i=1}^{j-1}P_i(x_i-x_{i+1}) + P_ja_j,\quad 1\le j\le k\cr
\phantom{aaaaaaaaaaaaaa}=\sum_{i=1}^{j-1}P_i(x_i-x_{i+1}) + \tilde P_{j-1}\overline a,\quad k+1\le j
\le n,\cr
=0,\quad j= n+1.\phantom{nnnnnnnn}\cr
}
$$
 Hence:
$$
\sum_{i=1}^kq_ix_i\ge \sum_{i=1}^kq_iy_i,\quad 1\le k\le n; \qquad \sum_{i=1}^{n+1}q_ix_i= \sum_{i=1}^{n+1}q_iy_i,
$$ 
and by  HLPKF  if $f$ is convex then
$$
\sum_{i=1}^{n+1}q_if(x_i)\ge  \sum_{i=1}^{n+1}q_if(y_i)=0,
$$
which is just $(J_n)$.
\QED
A variant of this result can be found in \ref{ABMP}.

In note on  Proof (\rn2)  that we do not use the full force of (P) as only the first and last weights are required to dominate the negative weights. While (P) makes much more demands on the negative weights  than does (S)  its real advantage, as Pe\v cari\'c  pointed out, is that no requirement of monotonicity of the elements of the \tup{n} is needed. This allowed an extension of Theorem 7 to convex functions of several variables as  we shall now demonstrate; \ref{MP}.

If $U\subseteq\Bbb R^k,\, k\ge 1$ where $U$ is a convex set then  the definition of convexity is, with a slight change in notation, just that given in (1): for all $\bf x,\bf y\in U$
$$
D(t)= D_2(t)= f\bigl((1-t){\bf x} + t{\bf y}\bigr)-\bigl((1-t)f({\bf x})+tf({\bf y})\bigr)\le 0,\quad 0\le t\le 1,
$$
and the convexity of $U$ ensures that $(1-t){\bf x} + t{\bf y}\in U$. Further one of the standard proof of $(J_n)$ can be applied in this situation to obtain Jensen's inequality for such functions $f$. Of course we cannot hope to extend the Steffensen result, if $k\ge 2$,  as the concept of increasing order of the points in $U$ is not available but the Pe\v cari\'c argument can be extended using the same proof as the one given above on the case $k=1$ and uses the same notations.

\proclaimsc{Theorem 8}{{\sl Let $U$ be an open convex set in $\BR^k$,
 ${\bf a}_i\in U, 1\le i\le n$, and  let $ p_i,  1\le i\le n$, be non-zero real  numbers with $P_n=1$ and  $I_- = \{i; 1\le i\le n\land p_i<0\}, I_+= \{i; 1\le i\le n\land p_i>0\}$. Further assume that  $\forall i,\; i\in I_-,\; \bf a_i$ lies in the convex hull of the set $\{{\bf a}_i; i\in I_+\}$ and that $\forall j,\; j\in I_+,\; p_j + \sum_{i\in I_-}\!p_i\ge0$.  If $f\colon U\mapsto \BR$ is convex then $(J_n)$ holds.}}\par
\qed Proof (\rn2) of Theorem 7 can be applied with almost no change although the notation ie a little messier.

If $i\in I_-$ then  for some $t_j^{(i)},\, 0\le t_j^{(i)}\le 1, \;  \sum_{j\in I_+}t_j^{(i)}=1$, ${\bf a}_i = \sum_{j\in I_+}t_j^{(i)}{\bf a}_j $ and so
$$
\eqalign
{
\overline {\bf a}=&\sum_{i=1}^np_i{\bf a}_i=\sum_{j\in I_+}p_j{\bf a}_j+\sum_{i\in I_-}p_i\Bigl(\sum_{j\in I_+}t_j^{(i)}{\bf a}_j\Bigr)\cr
=&\sum_{j\in I_+}\Bigl(p_j+\sum_{i\in I_-}p_it_j^{(i)}\Bigr){\bf a}_j\cr
=&\sum_{j\in I_+}q_j{\bf a}_j.\cr
}
$$
 where. as  in proof (\rn2) above, $0<q_j<1,\,\sum_{j\in I_+}q_j=1 $.  In this proof  we now use the full force of (P) and incidentally provide a needed proof that $\overline {\bf a}\in U$. The rest of the proof  proceeds as in proof (\rn2).\QED
Note that in the case $k=1$ the hypotheses imply that the smallest and largest element in the \tup{n}  have positive weights  each of which  dominates the sum of all the negative weights.

We now turn to $(\sim J)$ and note that proof (\rn4) of Theorem 7 can with a suitable change of hypotheses lead to this  inequality; \ref{P81a, P84a}.

\proclaimsc{Theorem 9}{\sl Let $ n, I, $ be as in Theorem 9 $p_1, \ldots p_n$ a real \tup{n}
\ with
$P_n=1$, then the reverse Jensen inequality holds for all functions $f$ convex on $I$  and for every monotonic with terms in $I$ if
and only if for some $m, 1\le m\le n, P_k\le 0, 1\le k< m$, and $\tilde P_k\le 0, m< k\le n$.}\par
\qed Looking at proof (\rn 4) of Theorem 7  we see that the present  hypotheses imply that 
$$
\sum_{i=1}^kq_ix_i\le \sum_{i=1}^kq_iy_i,\quad 1\le k\le n; \qquad \sum_{i=1}^{n+1}q_ix_i= \sum_{i=1}^{n+1}q_iy_i,
$$ 
and by  HLPKF  if $f$ is convex then
$$
\sum_{i=1}^{n+1}q_if(x_i)\le  \sum_{i=1}^{n+1}q_if(y_i)=0,
$$
which is just $(\sim J_n)$.\QED

 \sectn{Applications, Cases of Equality, Integral Results}
 The most obvious application so these extensions and reversals of the Jensen inequality are to mean inequalities. A large variety of means derive from the convexity of a particular function and  so we find tht these inequalities will now hold with negative weights under the above condition or are reversed. 
\subssn{An Example} If $ p_1, p_2, p_3, p_4$ are non-zero real numbers with $P_4=1$ and $a_1, a_2, a_3, a_4$ are distinct poitive numbers then, using the convexityv of the negative of the logarithmic function,  the particular case of (GA)
$$
a_1^{p_1}a_2^{p_2} a_3^{p_3} a_4^{p_4}\le p_1a_1 + p_2a_2 + p_3a_3 + p_4a_4 
$$
can be deduced  from Theorem 7 provided one of the following holds: 

\ \phantom{qq}(i) all the weights are positive;

\ \phantom{q}(ii) $a_1< a_2< a_3< a_4$ or $a_1> a_2>a_3> a_4$ and $0<p_1<1,0< P_2< 1, 0< P_3< 1$; 

\ \phantom{i}(iii) $a_1< a_2, a_3< a_4$, 
and $p_1> 0, p_4> 0$ and $P_3> 0, \tilde P_3> 0$.

 The reverse inequality 
$$
a_1^{p_1}a_2^{p_2} a_3^{p_3} a_4^{p_4}\ge p_1a_1 + p_2a_2 + p_3a_3 + p_4a_4, 
$$
can be deduced from Theorem 5 or Theorem 9  if one of the
following holds:

\ \phantom{qq}(i) only one of the weights is positive;

\ \phantom{q}(ii) either $a_1> a_2> a_3> a_4$, or $a_1<a_2, a_3< a_4$ and  either $0<p_1<1 0$ and $\tilde P_2,\tilde P_3, p_4< 0$, or  $0<p_2<1 0$ and $p_1,\tilde P_3< 0, p_4<0$, or $0<p_3<1 $ and $p_1,P_2, p_4< 0$ or $0<p_4<1 $ and $p_1, P_2, P_3<0$.
\subssn{The Pseudo Means of Alzer}
 A particular case of Theorem 5 has been studied by Alzer  under the name of pseudo-means, \refb{B03 pp. 171--173}, \ref{Al}.
\proclaimsc{Corollary 10}{ If $f$ is convex on $I$  and $p_i, 1\le i\le n$, are positive weights with $P_n=1$ then 
$$
f\Bigl(\recip{p_1}\bigl(a_1-\sum_{i=2}^np_ia_i\bigr)\Bigr)\ge \recip{p_1}\bigl(f(a_1)-\sum_{i=2}^np_if(a_i)\bigr),
$$
 provided $a_i,1\le i\le n,\drecip{p_1}\bigl(a_1-\sum_{i=2}^np_1a_1\bigr)\in I$.}\par
 A particular case when $ f(x) = x^{s/r},\, 0<r<s, \, x>0, $ leads to the inequality
$$
\Bigl(\recip{p_1}\bigl(a_1^s-\sum_{i=2}^np_ia_i^s\bigr)\Bigr)^{1/s}\ge \Bigl(\recip{p_1}\bigl(a_1^r-\sum_{i=2}^np_ia_i^r\bigr)\Bigr)^{1/r}.
$$
 A related topic is the Acz\'el-Lorenz inequalities; \refb{B03 pp. 198--199, PPT pp.124--126}, \ref{Ac}.
\subssn{The Inverse Means of Nanjundiah} Nanjundiah devised some very ingenious arguments using his idea of inverse means, \refb{B03 pp.136--137, 226}, \ref{N}. In the case of $ r>0$ Nanjundiah's inverse $r$-th power mean of order $n$ is defined as follows:  let $\bf a, w,$  be two sequences of positive numbers then
$$
\eufm N_n^{[r]}({\bf a; w})=\biggl(\frac{W_n}{w_n}a_n^r-\frac{W_{n-1}}{w_n}a_{n-1}^r\biggr)^{1/r}.
$$
An immediate consequence of Theorem 2 with $  f(x) = x^{s/r},\, 0<r<s, \, x>0,$ is the inequality
$$
\eufm N_n^{[r]}({\bf a; w})\ge \eufm N_n^{[s]}({\bf a; w}).
$$

\subssn{Comparable  Means}
If $\phi$ is a strictly increasing function then  a quasi-arithmetic mean is defined as follows:

$$ \eufm M_{\phi}( {\bf a;w})= \phi^{-1}\Bigl(\recip{W_n}\sum_{i=1}^nw_i \phi(a_i)\Bigr).
$$

An important question is when two such means are comparable, that is: when is it always true that: 
$$
\eufm M_{\phi}( {\bf a; w})\le \eufm M_{ \psi}( {\bf a; w})
$$

 Writing $ \phi(a_i)= b_i, 1\le i\le n,$ this last inequality :
$$
\psi\circ\phi^{-1}\Bigl(\recip{W_n}\sum_{i=1}^nw_i b_i\Bigr)\le\recip{W_n} \sum_{i=1}^nw_i \psi\circ\phi^{-1}(b_i),
$$
showing, from $(J_n)$, that the means are comparable exactly when $\psi\circ\phi^{-1}$  is convex, \refb{B03 pp. 273--277}.  Using Theorem 7 we can now allow negative weights in the comparison  and by using Theorem 5 or 9  get the opposite comparison; \ref{ABMP}.

Dar\'oczy \& P\'ales, \ref{DP},   have defined a  class of general means that they called {\sl L-conjugate means}:
$$
L_{\phi}^{{\eufm M}_1,\ldots {\eufm M}_n}({\bf a; u;v})=L_{\phi}({\bf a; u;v})
=\inv \phi\biggl(\sum_{i=1}^m{}u_i\phi(a_i)  - \sum_{j=1}^nv_j\phi\circ\eufm M_j({\bf a})\biggr)
$$
where $U_m-V_n= 1,\;  u_i>0,\, 1\le i\le m,\; v_j>0,\, 1\le j\le n$, $\eufm M_j, 1\le j\le n$ are means on \tups{m}  and $ \phi$ is as above.

Now suppose we wish to compare two L-conjugate means:
$$
L_{\phi}({\bf a; u;v})\le L_{\psi}({\bf a; u;v} ),
$$ 
Using the above substitution, $ \phi(a_i)= b_i, 1\le i\le m,$  and writing $\eufm N _j= \phi\circ \eufm M_j$ this last inequality becomes
$$
\psi\circ\inv\phi\Bigl(\sum_{i=1}^m{}u_ib_i-\sum_{j=1}^nv_j\eufm N_j({\bf b})\Bigr)\le \sum_{i=1}^mu_i\psi\circ\inv\phi(b_i)  -  \sum_{j=1}^nv_j\psi\circ\inv\phi\circ \eufm N_j({\bf b})
$$
 which, from  Theorem 8 in the case $k=1$, holds if $\psi\circ\inv\phi$ is convex, as for the quasi-arithmetic means; \ref{MP}.

In this sense this result  Pe\v cari\'c result gives a property of convex functions analogous to that of Jensen's inequality but useful for these means whereas Jensen's inequality is useful for the classical quasi-arithmetic  means.

It should be remarked that extensions of this comparison result can be obtained allowing the weights $\bf u, v$ to be real and using Theorem 7;   see \ref{ABPM}.

\subssn{Cases of Equality} Clearly the function $D$ of  (1)  is zero if either  $t=0$, $t=1$ or $x=y$; if otherwise $D<0$ then $f$ is said to be strictly convex. If this is the case then Jensen's inequality, $(J_n)$, is strict unless $a_1=\cdots =a_n$.

It follows easily from the proof of Theorem 5  that $(\sim J_n)$ holds strictly for strictly convex functions under the conditions of that theorem  unless $a_1=\cdots =a_n$.

In Theorem 7, Steffensen's extension of Jensen's inequality, the same is true by a consideration of proof (\rn2); see \ref{ABMP}.
\subssn{Integral Results} Most if not all of the above results have integral analogues but  a discussion of these would take us beyond the bounds of this paper; \refb{B03 p.371, PPT pp. 45--47, 84--87}, \ref{P81b}.
\sectn{Bibliography} 
\nr{[ABMP] S.\ Abramovich, M.\ Klari\v ci\'c Bakula, M.\ Mati\'c and J.\ Pe\v cari\'c} A variant of Jensen-Steffensen's inequality and quasi-arithmetic means, {J.\ Math.\ Anal.\ Appl.\/}, 307 (2005), 370--386.
\nr{[Ac] J.\ Acz\'el} Some general methods in the theory of functional equations in one variable. New applications of functional equations, {\sl Uspehi Mat.\ Nauk(N.S.)\/}, 11 (1956), 3--68
\nr{[Al] Horst Alzer}  Inequalities  for pseudo arithmetic and geometric means, {\sl General Inequalities\ Volume 6, Proceedings of the Sixth  International Conferences on General Inequalities\/}, 5--16,
Oberwolfach 1990.
\nr {[B98] P.\ S. Bullen} The Jensen-Steffensen inequality,  {\sl Math.\ Ineq.\ App.\/},1 (1998)m, 391--401. 
 \nr {[B03] P.\ S. Bullen}  {\sl Handbook of Means and Their Inequalities}, Kluwer Academic Publishers, Dordrecht, 2003.
\nr{[BB] Edwin F.\ Beckenbach and Richard Bellman} {\sl Inequalities\/},  Springer-Verlag, Berlin-Heidelberg-New
York, 1961.
\nr{[DP] Zolt\'an Dar\'oczy \& Zsolt P\'ales}  On a class of means of several variables,  {\sl\ Math.\ Ineq.\ App.\/}, 4 (2001),  
(20), 331--334.
\nr{[HLP] G.\ H.\ Hardy, J.\ E.\ Littlewood and G.\ P\'olya} {\sl  Inequalities\/},
Cambridge University Press, Cambridge,  1934.\ 
\nr{[LP] P.\ T.\ Landsberg \& J.\ E.\ Pe\v cari\'c} Thermodynamics, inequalities and negative heat capacities, {\sl\ Phys.\ Rev.}, A35 (1987),
 4397--4403.
\nr{[MO] Albert W.\ Marshall and  Ingram  Olkin} {\sl  Inequalities: Theory of Majorization and
Its Applications\/}, Academic Press,\  New York,\ 1979.
\nr{[MP] Anita  Matkovi\'c and Josip  Pe\v cari\'c} \ A variant of Jensen's inequality for convex functions of several variables, 
 {\sl J.\ Math.\ Ineq.\/}, 1 (2007), 45--51.
\nr{[MPF] D.\ S.\ Mitrinovi\'c, J.\ E.\ Pe\v cari\'c and  A.\ M.\ Fink} {\sl  Classical and New
Inequalities  in Analysis},  D\ Reidel, Dordrecht,   1993.
\nr{[N] T.\ S.\ Nanjundiah} Sharpening some classical inequalities, {\sl\ Math.\ Student}, 20 (1952), 24--25.
\nr{[P81a] J.\ E.\ Pe\v cari\'c} Inverse of Jensen-Steffensen inequality, {\sl\ Glasnik Mat.\/},16 (1981), 229--233. 
\nr{[P81b] Josip E.\ Pe\v cari\'c} A new proof of the Jensen-Steffensen inequality, {\sl\ Mathematica\  \  R\'ev.\ Anal.\ Num.\ Th\'eorie Appr.}, 23 (1981), 73--77. 
\nr{[P84a] J.\ E.\ Pe\v cari\'c} Inverse of Jensen-Steffensen inequality. (II), {\sl\ Glasnik Mat.\/},19 (1984), 235--238. 
\nr{[P84b] Josip  E.\ Pe\v cari\'c} A simple proof of the Jensen-Steffensen inequality, {\sl  Amer.\ Math.\ Monthly\/}, 91 (1984), 195--196.
\nr{[P90] Josip E.\ Pe\v cari\'c} Notes on Jensen's inequality, {\sl  General inequalities, 6 (Oberwolfach, 1990)},  449--454; Internat. Ser. Numer. Math., 103, BirkhŠuser, Basel, 1992 .
 \nr{[PPT] Josip E.\ Pe\v cari\'c, Frank Proschan and Y.\ L.\ Tong}, {\sl Convex Functions, Partial Orderings and Statistical Applications}, Academic Press, New York 1992.
 \nr{[RV] A.\ Wayne Roberts \& Dale E.\ Varberg} {\sl Convex Functions\/},   Academic Press, New
York-London,  1973.
\nr{[S] J.\ F.\ Steffensen} On certain inequalities and methods of approximation, {\sl  J.\ Inst.\
Actuar.\/}, 51 (1919), 274--297.
\nr {[VP] Petar M.\ Vasi\'c and Josip E.\ Pe\v cari\'c}  On the Jensen inequality, {\sl  Univ.\ Beograd Publ.\
Elektrotehn.\ Fak.\ Ser.\ Mat.\ Fiz.\/},  No.634--No.677, (1979), 50--54
\bigskip

\rightline{
 P.S.Bullen}

\rightline{Department of Mathematics}

\rightline{University of British Columbia}

\rightline{Vancouver BC} 

\rightline{Canada   V6T 1Z2}

\end